\documentclass{article}
\usepackage{graphicx} %

\usepackage[unicode,hyperfootnotes=false,hidelinks,hyperindex=false]{hyperref}
\usepackage{graphicx,tikz,framed,amssymb,multicol,enumerate,xcolor}
\usepackage{amsmath,amsthm,amsfonts,amssymb,xcolor,etoolbox,cancel}
\usepackage[english]{babel}
\usepackage[numbers]{natbib}
\usepackage[T1]{fontenc}
\usepackage{dsfont}
\graphicspath{ {./logos/} }
\usepackage[export]{adjustbox}

\usepackage[symbol]{footmisc}
\usepackage{perpage}
\MakePerPage{footnote} %

\numberwithin{equation}{section}

\usepackage[arrow,matrix]{xy}
\newlength{\xywd}
\newcommand{\xymapsto}[2][]{\mathrel{%
  \sbox{0}{$\scriptstyle#1$}%
  \xywd=\wd0
  \sbox{0}{$\scriptstyle#2$}%
  \ifdim\wd0>\xywd \xywd=\wd0 \fi
  \xymatrix@M=0ex@C\dimexpr\xywd+0.7em\relax{{}\ar@{|->}[r]^{#2}_{#1}&{}}%
}}
\newcommand{\xyrightarrow}[2][]{\mathrel{%
  \sbox{0}{$\scriptstyle#1$}%
  \xywd=\wd0
  \sbox{0}{$\scriptstyle#2$}%
  \ifdim\wd0>\xywd \xywd=\wd0 \fi
  \xymatrix@M=0ex@C\dimexpr\xywd+0.7em\relax{{}\ar@{->}[r]^{#2}_{#1}&{}}%
}}

\newcommand{\N}{\mathbb{N}}

\newcommand{\defi}{\mathrel{\mathop:}=}

\newcommand{\union}{\mathop{\textstyle\bigcup}}

\DeclareMathOperator{\ran}{ran}

\newcommand{\WO}{\mathrm{WO}}
\newcommand{\estf}{\mathrm{EST}}
\newcommand{\zfc}{\mathrm{ZFC}}
\newcommand{\esti}{\text{EST+WAI}}

\newcommand{\axiomas}[1]{\mathrm{#1}}

\newcommand{\AC}{\axiomas{AC}}

\newcommand{\ZFC}{\axiomas{ZFC}}
\newcommand{\ZF}{\axiomas{ZF}}

\newtheorem{theorem}{Theorem}[section]
\newtheorem{lemma}[theorem]{Lemma}

\newtheorem{prop}[theorem]{Proposition}
\newtheorem{corollary}[theorem]{Corollary}

\newtheorem*{WeakInf}{Weak Axiom of Infinity}

\theoremstyle{definition}
\newtheorem{definition}[theorem]{Definition}

\theoremstyle{remark}

\newtheorem*{claim*}{Claim}

\newtheorem*{remark*}{Remark}

\newenvironment{hint*}{\pushQED{\qed}\begin{ayuditaaa*}}{\popQED\end{ayuditaaa*}}
\newtheorem*{solveitnow*}{Solución} %
{%
  \pushQED{\qed}\begin{solveitnow*}}
{\popQED\end{solveitnow*}}

\usepackage{tocloft}
\newlength{\imarg}

\tikzset{
every picture/.append style={
  execute at begin picture={\deactivatequoting},
  execute at end picture={\activatequoting}
  }
}

\hypersetup{
  pdftitle={Hartogs},
  pdfauthor={Pedro Sanchez Terraf, Joel Kuperman, Gervasio Figueroa},
  pdfsubject={Set Theory},
  pdfkeywords={axiomas, funciones, límites, sucesiones, derivadas},
  colorlinks,
  linkcolor={blue!35!black},
  citecolor={blue!35!black},
  urlcolor={blue!35!black}
}

\title{Iterating logarithms and well-foundedness}
\author{
  Gervasio Figueroa\thanks{
    Universidad Nacional de Córdoba. 
    Facultad de Matemática, Astronomía,  Física y
    Computación. \\
    \indent \phantom{*} Centro de Investigación y Estudios de Matemática (CIEM-FaMAF).
    Córdoba. Argentina.\\
    Supported by Secyt-UNC project 33620230100751CB.
  }
  \and
  Joel Kuperman\footnotemark[1]
  \and
  Pedro Sánchez Terraf\footnotemark[1]
}

\begin{document}
\maketitle

\begin{abstract}
  We prove that there is no sequence $\{ X_n \}_{n\in\omega}$ such that for all
  $n$, there is an injection from $\mathcal{P}(X_{n+1})$ to $X_n$ under a
  restricted set of axioms that conforms to the “working mathematician's”
  intuition (excluding Choice, Replacement, and Foundation).

  We also show such a sequence cannot be defined using a formula either. More
  precisely, it is inconsistent with Zermelo set theory minus Infinity to have a
  formula $\Phi(\cdot,\cdot)$ such that for each $n\in\omega$ there are $X$ and
  $Y$ with $\Phi(n,X)$, $\Phi(n+1,Y)$, and  $\mathcal{P}(Y)$ being injectable
  into $X$.

  We finally show that the theorem fails if we allow infinitely many
  definitions, one for each term of the sequence, by using Hamkins, Linetsky,
  and Reitz's pointwise definable extensions of countable models of $\ZFC$.
\end{abstract}

\section{Introduction}
Well-foundedness undoubtedly plays a decisive role throughout Set Theory. This
paper intends to be a contribution to the question of how much of this notion persists
when considering restricted sets of axioms that conform to the “working
mathematician's” intuition. At the same time, we are interested in examples of
elementary (or easy-to-state) problems, solvable by the use of set-theoretic
axioms that do not defy any conservative intuition, but require notions beyond
naive concepts.

Take, for example, the notion of cardinality. It is equivalent to $\AC$ over
$\ZF$ that cardinalities are well-\emph{ordered}; notably, it is an open
question \cite[p.~312]{zbMATH07977883} \cite[Prop.~4.2]{2023arXiv231009643B} whether the fact that the universe of sets is well-founded under the
relation defined by
\begin{quote}
  $A \leq_c B$ if and only if there is an injection from $A$ to $B$,
\end{quote}
implies $\AC$.

Nevertheless, we consider the Axiom of Choice to lie beyond our list of
“intuitive” axioms, which will also exclude Foundation. In this restricted
setting, it is clear that we cannot expect well-foundedness of cardinalities but
we still have that cardinalities cannot decrease at an iterated-exponential rate
(much in the same way that it is impossible to iterate $\log_2$ on the real
numbers indefinitely). Namely, $\ZF$ proves \cite[Exr.~I(12)]{kunen1980} that
\begin{theorem}\label{informal_hartogs_set}
  There is no sequence $\{X_n\}_{n\in\omega}$ such that  for all $n$, there is an injection
    from $\mathcal{P}(X_{n+1})$ to $X_n$.
\end{theorem}

In this article we restrict ourselves a little more and prove the above fact
(Theorem~\ref{hartogs_set}) without using the Replacement axiom, motivated by the discussion in \cite{PedroStack}.

Moreover, we show (Theorem~\ref{hartogs_class}) that not only this sequence
cannot exist as a set, it cannot even be defined by a formula (this is relevant
due to the absence of Replacement).
More precisely, it is inconsistent with Zermelo set theory minus Infinity to have a formula $\Phi(\cdot,\cdot)$ such that for each $n\in\omega$ there are $X$ and $Y$ with
\[
	\Phi(n,X) \land \Phi(n+1,Y) \land \mathcal{P}(Y) \leq_c X.
\]

Finally, we prove (Theorem~\ref{th:nonuniform-exp-decreasing}) that the theorem fails if we allow infinitely many
definitions, one for each term of the sequence. That is, for
each $n$ there is a formula $\Phi_n$ that defines a set and such that the set of
formulas given by
\[
  \forall X\, \forall Y,\,  \Phi_n(X) \land \Phi_{n+1}(Y) \implies
  \mathcal{P}(Y) <_c X
\]
for each $n$, is consistent with $\zfc$.

The proof of the last theorem required us to extend a non well-founded countable
model of $\zfc$ to a pointwise definable one by forcing, which \cite{hamkins}
shows it can always be done.

\section{Preliminaries}
We begin by laying out the list of axioms that we believe represents one possible
description of the intuitive concept of set.
\begin{definition}
	We shall use $\estf$ (Elementary Set Theory) to mean the axioms of Extensionality, Comprehension, Pairing, Union and Powerset.
\end{definition}
Since we are leaving
Replacement aside, one other axiom that might be considered affected is the
Axiom of Infinity. For instance, it is a well known fact that Zermelo's
\textit{Unendlichen} axiom does not imply Von Neumann's one without Replacement
(e.g., \cite[Exercise II.4.21]{Kunen}), and the lack of $\AC$ brings even more
elements to the picture (infinite Dedekind-finite sets, amorphous sets, etc.).

We hence prefer to consider a version of this axiom that does not involve
ordinals in any implicit sense. Thus, we state a version of Infinity weaker than the
existence of $\omega$. Before stating this axiom, let us recall the following
definitions which we also use to fix some notations.
\begin{definition}
	\begin{enumerate}
		  \item For sets  $X$ and $Y$, $X\leq_c Y$ iff there is an injective function $f:X\to Y$, $X =_c Y$ iff there is a bijective function $f:X\to Y$, and $X <_c Y$ iff $X\leq_c Y$ and $X\neq_c Y$.
		\item A set $n$ is a \emph{natural number} iff it is the empty set or a successor ordinal and for each $m\in n$, $m$ is the empty set or a successor ordinal.
		\item A set $A$ is \emph{finite} iff $A\leq_c n$ for some natural number $n$.
		\item A set is \emph{infinite} iff it is not finite.
	\end{enumerate}
\end{definition}
In the following, we write “$n\in\omega$” for the formula defining natural numbers:
\[
n \in \text{Ord} \land \forall x, \, x \in n \lor  x = n \implies x = \emptyset \lor (\exists y \in
n, \, x = y \cup \{y\}).
\]

\begin{WeakInf}
  \textup{(WAI)} There is an infinite set.
\end{WeakInf}
By using \cite[Exercise II.4.21]{Kunen}, we know that $\esti$ is consistent with
the negation of the standard Axiom of Infinity relative to $\zfc$.

\begin{definition}
  For linear orders $V$ and $U$,
  \begin{enumerate}

  \item
    $V\leq_o U$ iff  $V$ is
    isomorphic to an initial segment of $U$;
  \item
    $V<_o U$ when  $V\leq_o U$ and  $V$ is not order isomorphic
    to $U$.
  \end{enumerate}
\end{definition}

The following proposition follows from $\estf$ and is necessary to prove some elementary facts about finite sets using induction over $\omega$, which may be a proper class even in $\esti$.

\begin{prop}[{\cite[Corollary 7.33]{Moscho}}]\label{class_wo}
  Every non-empty class $\mathcal{C}$ of well ordered sets has a $\leq_o$-least member.
\end{prop}

We start by laying out a few well known definitions and propositions about the Hartogs function. Of particular note is the bound $\aleph(X) <_c \mathcal{P}^{3}(X)$, which is the key idea for our proof of Theorem \ref{informal_hartogs_set}.

\begin{definition}\label{aleph}
	Let $X$ be any set,
	\[
	\WO(X):=\{\langle A,R\rangle \mid A\subseteq X \ \text{is well-ordered by}\ R\}, 
	\]
	and $\aleph(X):= \WO(X)/{\cong}$ where $\cong$ denotes order isomorphism.
	Given $\aleph(X)$, and equivalence classes $A,B \in\aleph(X)$ we shall say $A
        \leq_{\aleph(X)} B$ iff for any (all) $U\in A$ and $V\in B$, $U\leq_o V$
        holds.
\end{definition}

From Moschovakis \cite[Thm.~7.34]{Moscho} we conclude
\begin{lemma}\label{aleph2}
	For all $X$, $\aleph(X)$ has the following properties:
	\begin{enumerate}
		\item $\aleph(X)\nleq_c X$.
		\item $(\aleph(X),\leq_{\aleph(X)})$ is a well order and is not
                  equinumerous with any of its proper initial segments.
		\item If $W$ is a well-orderable set and $W \nleq_c X$, then $\aleph(X) \leq_o W$.\qed
	\end{enumerate}
\end{lemma}

\begin{corollary}\label{cor:ltc-iff-lto}
  $\aleph(X) <_c \aleph(Y)$ if and only if $\aleph(X) <_o \aleph(Y)$.
\end{corollary}
\begin{proof}
  It follows by Lemma~\ref{aleph2} and the fact that the relation $\leq_o$ is
  total (Proposition~\ref{class_wo}, or more precisely by \cite[Thm.~7.31]{Moscho}).
\end{proof}

\begin{theorem}\label{thm:partes-cubo}
	For any set $X$, $\aleph(X) <_c \mathcal{P}^3(X)$.
\end{theorem}
\begin{proof}	
  Let $f:\WO(X)\rightarrow \mathcal{P}^2(X)$ defined by $f(\langle
  A,\mathrel{R}\rangle) = \{ \{ b\mid b \mathrel{R} a \} \mid a\in A \}$. $f$
  is clearly injective. Hence, the injective function
  $\Tilde{f}: \mathcal{P} (\WO(X)) \rightarrow \mathcal{P}^3(X)$ given
  by $\Tilde{f}(S) = \{ f(s) \mid s\in S \}$ witnesses
  $\mathcal{P} (\WO(X)) \leq_c \mathcal{P}^3(X)$. Since
  $\aleph (X) \subset    \mathcal{P} (\WO(X))$ we conclude
  $\aleph(X) \leq_c \mathcal{P}^3(X)$.
            
	Let us prove that $\aleph(X) \neq_c \mathcal{P}^3(X)$. 
	By way of contradiction, we assume they are equinumerous; then $\mathcal{P}^{3}(X)$ is well-orderable. As an immediate consequence we have that $\mathcal{P}^{2}(X)$ is well-orderable. But then $\mathcal{P}^{2}(X)$ is a well-orderable set and it satisfies 	
	$\mathcal{P}^{2}(X)\nleq_c X$.
	So Lemma \ref{aleph2} implies that 	
	$\aleph(X)\leq_c\mathcal{P}^{2}(X)$.
	But 	
	$\mathcal{P}^{3}(X) =_c \aleph(X)$	
	by assumption, so we have shown $\mathcal{P}^{3}(X) \leq_c \mathcal{P}^{2}(X)$, which gives us our contradiction.
\end{proof}

While the full strength of $\zfc$ proves $\forall X \ \aleph(X) <_c \mathcal{P}^{2}(X)$, the above theorem is the best bound we can get in $\ZF$: It is consistent with $\ZF$ that the inequality $\aleph(X) <_c \mathcal{P}^2(X)$ fails for some $X$. More specifically, for any ordinal $\alpha$ such that $\alpha < \aleph_{\alpha}$ it is consistent with $\ZF$ to have $\aleph(X) =_c \aleph_\alpha$ and for $\mathcal{P}^{2}(X)$ and $\aleph(X)$ to have incomparable cardinalities \cite[Theorems 10 and 11]{hickman};
see also the discussion at \cite{stack_powerset}.

Using the well-foundedness of $\leq_o$, we can verify that the principle of
ordinary Induction can be proved under $\estf$:
\begin{theorem}
	Let $P(n,\vec{x})$ be a first order formula in the language of set theory whose free variables are either $n$ or contained in $\vec{x}$. Assume that the following statements hold:
	\begin{enumerate}[(i)]
		\item $P(\emptyset,\vec{x})$.
		\item\label{inductive} Let $n$ be a natural number, then $P(n,\vec{x}) \implies P(n\cup \{n\},\vec{x})$.
	\end{enumerate}
	Then $P(n,\vec{x})$ is true for every natural number.
\end{theorem}
Having Induction, the following elementary results become available.
\begin{lemma}
	Let $X$ be an infinite set and $n$ a natural number, then $n <_c X$.
\end{lemma}

\begin{corollary}\label{cor:finite-smaller-infinite}
	Let $A$ be a finite set and $X$ an infinite set, then $A <_c X$.
\end{corollary}

Though we are working in an axiomatic system too weak to prove the existence of
the ordinal $\omega$, nevertheless it is handy to have a first order definition
of a similar well order.
\begin{definition}
  A nonempty well order has \emph{type $\omega$} iff  every element
  has an immediate successor, and every element (except the first one) has an
  immediate predecessor.
\end{definition}
It can be seen in $\estf$ that this actually defines an order type (i.e.,
having type $\omega$ is an order isomorphism class) and coincides with being
isomorphic to $\omega$ in case it exists.

The existence of an order of type $\omega$ might look stronger than the Weak
Axiom of Infinity. It can actually be seen under $\estf$ that they are equivalent
by using Hartogs' Lemma.

In the following lemma, $\text{Succ}(x,y)$ denotes that $y$ is an immediate
successor of $x$ under $<$: 
\[
  x<y \land \forall a\in A,\, a < y \implies (a=x \lor a < x).
\] 
\begin{lemma}\label{lem:omega-segment}
  If $\langle X,< \rangle$ is an infinite well order that does not have type
  $\omega$, then there is some $x\in X$ such that the restriction of $<$ to
  $x{\downarrow}$ (the initial segment determined by $x$) has type $\omega$.
\end{lemma}
\begin{proof}
	Since $\langle X,< \rangle$ does not have type $\omega$, there must be some $x\in X$ that is not the first element of $X$ and satisfies
	\[
	\neg(\exists y\in X,\, \text{Succ}(y,x)) \lor \neg (\exists z\in X,\, \text{Succ}(x,z)).
	\]
	We shall see that in both cases $x{\downarrow}$ is infinite. Let us assume 
	\[
	\neg\exists y\in X,\, \text{Succ}(y,x),
	\]
	then $x{\downarrow}$ must be infinite. It cannot be empty, because the first element of $X$ belongs to $x{\downarrow}$; and if it were finite, it would have a maximum $y$, which would imply $\text{Succ}(y,x)$.
	
	On the other hand, if we assume
	\[
	\neg \exists z\in X,\, \text{Succ}(x,z),
	\]
	then $X = x{\downarrow} \cup \{x\}$. Given that $X$ is infinite by assumption, it follows that $x{\downarrow}$ must be infinite.
	
	Finally, let us take $Y = \{x\in X \mid x{\downarrow} \text{ is infinite}\}$. This set is non empty by the previous argument, so it has a least element $y$. We will prove that $y{\downarrow}$ has type $\omega$ with the restricted order.

	Let $a\in y{\downarrow}$; since $y$ is the first element whose set of
        predecessors is infinite, $a{\downarrow}$ is necessarily finite. In
        particular it is either empty, in which case $a$ is the first element of
        $y{\downarrow}$, or it has a greatest element $b$ in which case we have $\text{Succ}(b,a)$. All that remains to verify is the existence of some $c\in y{\downarrow}$ such that Succ$(a,c)$.
	First let us note that $a<y$ and $\neg\text{Succ}(a,y)$, so the set $\{z\in y{\downarrow} \mid a<z\}$ is non empty. Taking its first element we obtain a $c$ with the desired property.
\end{proof}

\begin{prop}[$\esti$]\label{pr:well-order-type-omega}
	There is a well ordering of type $\omega$.
\end{prop}
\begin{proof}
	Let $X$ be any infinite set. By the previous
        lemma, it is enough to show that $\aleph(X)$ is not finite. If $\aleph(X)$ were finite, then we would
        have $\aleph(X) <_c X$ by Corollary~\ref{cor:finite-smaller-infinite},
        which contradicts Lemma~\ref{aleph2}.
\end{proof}

In order to state our main first theorem, we use Peano systems as a substitute for $\omega$.
\begin{definition}[{\cite[Definition 5.1]{Moscho}}]\label{peano}
	A \emph{Peano system} is any structure $(\mathcal{N},o,S)$ which satisfies the following properties:
	\begin{enumerate}[(1)]
		\item $\mathcal{N}$ is a set (the \emph{domain} of the Peano system), and $o\in\mathcal{N}$.
		\item $S$ is a function on $\mathcal{N}$, $S: \mathcal{N}\rightarrow\mathcal{N}$.
		\item $S$ is an injection.
		\item For each $n\in\mathcal{N}$, $S(n) \neq o$.
		\item For each $X\subseteq\mathcal{N}$,
		\[
		  o\in X \land (\forall n\in \mathcal{N}, \, n\in X \implies S(n)\in X) \implies X = \mathcal{N}.
		\]
	\end{enumerate}
\end{definition}

Note that the existence of a Peano system immediately implies the Weak Axiom of Infinity.
Assuming there is a Peano system, $\estf$ proves that one can define functions recursively (\cite[Theorem 5.6]{Moscho}). From there, it is straightforward to define addition in a given Peano system as $n+o = n$ and $n+S(m) = S(n+m)$ for each $n,m$.

The following lemma highlights the relation between Peano systems, well orderings, and the Weak Axiom of Infinity.

\begin{lemma}
	\label{peano_is_order}
	The following statements are equivalent:
	\begin{enumerate}[(1)]
		\item There is an infinite set (WAI).
		\item There is a well order of type $\omega$.
		\item There is a Peano system.
	\end{enumerate}
\end{lemma}
\begin{proof}
	Note that $(1){\implies}(2)$ is just
        Proposition~\ref{pr:well-order-type-omega}, and the converse is immediate.
	Let us prove $(2) {\implies} (3)$.

        Assume $\langle A,< \rangle$ is a
        well order of type $\omega$ and let $a_0$ be its first element.
        We can use Comprehension on the set $A\times A$: we simply take $S_A$ as $\{(x,y)\in A\times A: \text{Succ}(x,y)\}$. Let us now show that $(A,a_0,S_A)$ is indeed a Peano system. Properties (1) through (4) of Definition \ref{peano} are immediate consequences of the hypothesis on $\langle A,< \rangle$. Let $X\subseteq A$ be a set such that 
        \begin{equation}\label{eq:X-inductive}
          a_0\in X \land (\forall n\in A, \, n\in X \implies S_A(n)\in X).
	\end{equation}
	Assume $X\neq A$ and let $a'$ be the first element of $A\setminus X$. Using that $A$ has type $\omega$, we have that $a'$ is the first element of $A$ or there is some $b\in A$ such that $a' = S_A (b)$. In the first case we would have $a' = a_0$, which contradicts (\ref{eq:X-inductive}). So there must be some $b\in A$ such that $S_A (b) = a'$. But $a'$ is the first element of $A\setminus X$ and $b<a'$, which implies $b\in X$. Finally $S_A(b) \in X$, a contradiction which arose from assuming $A\setminus X$ was non-empty.
	
	To show $(3){\implies}(2)$, let $(\mathcal{N}, o, S)$ be a Peano system and define $n<m$ iff $\exists k\in\mathcal{N}{\setminus} \{o\}$ such that $n+k = m$.
	 $\langle \mathcal{N},< \rangle$ is a well ordering by \cite[Theorem 5.23]{Moscho}. The fact that it also has type $\omega$ is trivial from the definition of $<$.
\end{proof}

\section{Exponentially decreasing sequences of sets}

In this section we state and prove in $\estf$ a version of
Theorem~\ref{informal_hartogs_set} for sequences indexed by any Peano
system (as we have seen, the existence of such a system follows from WAI). In this way, we morally dispense with the Replacement instances needed to construct
$\omega$. To avoid cluttering the notation with the successor function $S$ of a
Peano system, we shall refer to $S(n)$ as $n+1$, $S(S(n))$ as $n+2$ and so on.

\begin{theorem}[$\estf$]
	\label{hartogs_set}
	Let $(\mathcal{N},o,S)$ be a Peano system. There is no sequence $\{X_n\}_{n\in \mathcal{N}}$ such that for all $n\in\mathcal{N}$, $\mathcal{P}(X_{n+1}) \leq_c X_n$.
\end{theorem}
\begin{proof}
	Assume there is such a sequence, 
	$\mathcal{S} := \{X_n\}_{n\in \mathcal{N}}$,
	then $\estf$ proves that its range $\{X_n \mid n\in \mathcal{N}\}$ exists.
        Let $\mathit{PWO} \defi \mathcal{P}(\WO(\union \{X_n \mid n\in \mathcal{N}\}))$.
	For all $X\in\ran(\mathcal{S})$, we have 
	\[
	\aleph(X) \subseteq \mathcal{P}(\WO(X)) \subseteq \mathit{PWO},
	\]
	where both inclusions follow from Definition \ref{aleph}.
	
	\noindent Thus, for all $X\in \{X_n \mid n\in \mathcal{N}\}$, we have
        $\aleph(X)\in\mathcal{P}(\mathit{PWO})$. Hence, $\esti$ proves the
        existence of the sequence $\{\aleph(X_n)\}_{n\in \mathcal{N}}$ by
        separating on $ \mathcal{N}\times\mathcal{P}(\mathit{PWO})$.

        Since applying powersets preserves lax cardinal inequalities, we have
        \[
          \mathcal{P}^3(X_{n+3}) \leq_c \mathcal{P}^2(X_{n+2})\leq_c \mathcal{P}(X_{n+1}) \leq_c X_{n}.
        \]
        By applying Theorem \ref{thm:partes-cubo}, we conclude
        $\aleph(X_{n+3}) <_c \mathcal{P}^{3}(X_{n+3}) \leq_c X_n$.
	
	Let us now show that $\aleph(X_{n+3}) <_c \aleph(X_n)$ for all $n\in
        \mathcal{N}$. Since both terms are $\leq_o$-comparable by Proposition~\ref{class_wo}, we have that either
        $\aleph(X_{n+3}) <_c \aleph(X_n)$ or $\aleph(X_n) \leq_c
        \aleph(X_{n+3})$ holds. But the latter statement implies $\aleph(X_n) <_c X_n$, which contradicts Lemma \ref{aleph2}.
	
	Finally, $\{\aleph(X_n) \mid n\in \mathcal{N}\}$
	is a class of well ordered sets that is not well founded with the $<_c$
        relation and hence not $\leq_o$-well founded, contradicting Proposition~\ref{class_wo}.
\end{proof}

The following theorem is similar to Theorem~\ref{hartogs_set}, but has a
somewhat different flavor, in that it works with the (possibly proper) class of
natural numbers and with a definable sequence of sets; actually, we allow
multiple options for each natural argument. In its statement and
proof,  $n+1$ denotes $n \cup \{n\}$.
\begin{theorem}[$\estf$]
	\label{hartogs_class}
	There is no formula $\Phi(\cdot,\cdot)$  with at least two free
        variables that satisfy both of the following properties:
	\begin{enumerate}[(i)]
		\item \label{nat_define} For each $n\in\omega$ there is a set $X$ that satisfies $\Phi(n,X)$.
		\item \label{dec_seq} For each $n\in\omega$ the following holds: 
		\[
			\Phi(n,X_0) \land \Phi(n+1,X) \implies \mathcal{P}(X) \leq_c X_0.
		\] 
	\end{enumerate}
\end{theorem}

\begin{proof}
	By way of contradiction, let $\Phi(\cdot,\cdot)$ be a formula with at
        least two free variables and assume it satisfies (\ref{nat_define}) and
        (\ref{dec_seq}), we shall take the following class (the range of the
        composition “$\aleph\circ \Phi$”):
	\[
		 \mathcal{C} := \{Y \mid  \exists n\in\omega,\, \exists X,\, \Phi(n,X) \land Y = \aleph(X)\}.
	\]
	Using (\ref{nat_define}), Definition \ref{aleph}, and Lemma \ref{aleph2}, $\mathcal{C}$ is a non-empty class of well ordered sets. By Proposition \ref{class_wo}, there must be a $\leq_o$-least member of $\mathcal{C}$.
	Let $Y_0$ be such a member of $\mathcal{C}$, then there must be some $X_0$ and $n_0 \in\omega$ such that $Y_0 = \aleph(X_0) \land \Phi(n_0,X_0)$.
	The assumption (\ref{nat_define}) implies there are $X_1 , X_2, X_3$ such that $\Phi(n_0 + i , X_i)$ for $i\in\{1,2,3\}$.
	The proof can be finished analogously to that of
        Theorem~\ref{hartogs_set} by applying (\ref{dec_seq}) for $X_i$ ($i\in\{0,1,2,3\}$).
\end{proof}

The previous theorem handles the case of a uniformly
definable sequence of sets. One possible next step in generality is that of a
(meta-)sequence of definable sets, each one given by a separate formula $\Phi_n$
with exactly one free variable, for each natural number $n$ in the
metatheory.

It turns out that a result analogous to Theorem~\ref{hartogs_class} in this
setting fails. To show this, we take advantage of
Hamkins, Linetsky, and Reitz's pointwise definable extensions to obtain an appropriate
meta-sequence of formulas.
\begin{theorem}[{\cite[Theorem~3(3)]{hamkins}}]\label{th:pointwise-extension}
  Every countable model of $\zfc$ has a class forcing extension that is
  pointwise definable.
\end{theorem}

Then our refutation takes this shape:
\begin{theorem}\label{th:nonuniform-exp-decreasing}
  Assume $\zfc$ is consistent. There is a collection of formulas
  $\Phi_n(\cdot)$ $(n\in\mathbb{N})$ with one free variable
  such that
  \begin{itemize}
  \item  for every $n$,
    $\estf \vdash \exists! X, \, \Phi_n (X)$; and
  \item The collection of all formulas $\forall X\, \forall Y,\,  \Phi_n(X) \land \Phi_{n+1}(Y) \implies
    \mathcal{P}(Y) <_c X$ with $n\in\mathbb{N}$ is consistent with $\zfc$.
  \end{itemize}
  In particular,  $\zfc \nvdash \exists X\,\exists Y,\, \Phi_n (X) \land \Phi_{n+1} (Y)
  \land \neg\bigl(\mathcal{P}(Y) <_c X\bigl)$ for each $n$.
\end{theorem}
\begin{proof}
  We start by taking a countable model of $\zfc$. 
  We extend its signature by adding countably many names for constants
  $c_n$, and for each $n\in\N$ we take the sentence $P_n$ to be
  \[
    c_{n+1} \text{ is a natural number} \land 2^{c_{n+1}} < c_n,
  \]
  where exponentiation of natural numbers is construed by its recursive
  definition and hence it is $\Delta^{\ZF}_1$.
  Since $\zfc+\{P_n\}_{n\in\N}$ is finitely satisfiable,
  it has a model $\mathcal{M}_c =
  (M,E,\{ x_n \}_{n\in \N})$ by Compactness (here, $x_n$ is the
  $\mathcal{M}_c$-interpretation of $c_n$). By
  Löwenheim-Skolem, we may assume that $\mathcal{M}_c$ is countable and
  let $\mathcal{M}=(M,E)$ be its reduct to the language of set
  theory. Finally, we extend it to a pointwise definable model $\mathcal{M}'=(M',E')$
  by using Theorem~\ref{th:pointwise-extension}.
  
  Let $\Omega_n$ be a definition of $x_n$ in $\mathcal{M}'$.	
  Then 
  \[
    \Phi_n (X) \defi \bigl((\exists! Y \ \Omega_n(Y))  \land \Omega_n(X)\bigr) \lor
    \bigl(\neg(\exists! Y \ \Omega_n(Y)) \land X=\emptyset\bigr)
  \] 
  satisfies that 
  \[
    \estf \vdash \exists! X \ \Phi_n (X).
  \]
  By construction, we have that for all $x\in M'$,
  \begin{equation}
    \label{eq:Phi-Omega}
    x = x_n  \quad\text{iff}\quad \mathcal{M}' \models \Phi_n[x]  \quad\text{iff}\quad \mathcal{M}'
    \models \Omega_n[x].
  \end{equation}
  Since the formula
  \[
    \varphi(X,Y) \defi \text{$X$ and $Y$ are natural numbers} \land 2^Y < X
  \]
  is $\Delta^{\ZF}_1$ and $\mathcal{M}'$ is a forcing extension of $\mathcal{M}$,
  $\varphi$ is absolute between them. Hence, for each $n\in \N$ the
  following sentence holds in $\mathcal{M}'$:
  \[
    \forall X\, \forall Y,\,  \Omega_n(X) \land \Omega_{n+1}(Y) \implies \varphi(X,Y).
  \]

  Finally, since $\zfc$ proves that $2^Y < X$ is equivalent to $\mathcal{P}(Y) <_c X$
  for $X,Y$ natural numbers, by (\ref{eq:Phi-Omega}) we have
  \[
    \mathcal{M}' \vDash \forall X\, \forall Y,\,  \Phi_n(X) \land \Phi_{n+1}(Y) \implies \mathcal{P}(Y) <_c X
  \]
  for each $n\in\N$, which shows the second bullet point.
\end{proof}

We have just seen that a non-uniform version of Theorem~\ref{hartogs_class}
fails. We recall that the proof of that theorem proceeded by obtaining a decreasing sequence of
well ordered sets from an exponentially decreasing one. Note that
Theorem~\ref{th:nonuniform-exp-decreasing} directly constructs a sequence
of the latter kind.

Apparently, this does not preclude the possibility of a model of
$\estf$ containing decreasing sequences of well orders but no exponentially
decreasing sequence of sets. Nevertheless, the counterexample constructed above
cannot distinguish between these two situations, since in any model having
an infinite decreasing sequence of natural numbers (and hence nonstandard
$\omega$), there exists an exponentially decreasing sequence of numbers
(given a nonstandard $c\in\omega$, take $d-1$ where $d\in\omega$ is the least
such that $c\leq 2^d$).  

\paragraph*{AI use disclosure}
AI was not used during the research nor the writing of this work.

\paragraph*{Acknowledgment}
We thank Prof.~Zachiri McKenzie (University of Chester) for discussions around
Proposition~\ref{pr:well-order-type-omega}, and Vicente Schkolnik (Universidad
Nacional de Córdoba) for the proof of strictness in
Theorem~\ref{thm:partes-cubo}. We also want to thank Nancy Moyano (Centro de
Investigación y Estudios de Matemática, CIEM-FaMAF) for support during this
project.

\bibliographystyle{plain}
\bibliography{bibliografia-hartogs}

@book{Kunen,
    title = {Set Theory},
    author = {Kenneth Kunen},
    isbn = {9781848900509},
    year = {2013},
    publisher = {College Publications}
    }

@book{Moscho,
	title = {Notes on Set Theory},
	author = {Yiannis Moschovakis},
	isbn = {978-0387-28722-5},
	year = {2006},
	publisher = {Springer}
}

@article{zbMATH07977883,
 author = {Tachtsis, Eleftherios},
 title = {No decreasing sequence of cardinals in the hierarchy of choice principles},
 fjournal = {Notre Dame Journal of Formal Logic},
 journal = {Notre Dame J. Formal Logic},
 issn = {0029-4527},
 volume = {65},
 number = {3},
 pages = {311--331},
 year = {2024},
 language = {English},
 doi = {10.1215/00294527-2024-0020},
 zbMATH = {7977883}
}

@MISC {PedroStack,
	TITLE = {{(ZF - Foundation)} proves that $\log_2$ can't be infinitely iterated: Alternative proof},
	AUTHOR = {Sánchez Terraf, Pedro},
	HOWPUBLISHED = {Mathematics Stack Exchange},
	NOTE = {\url{https://math.stackexchange.com/q/1422718} (version: 2017-04-13)},
	EPRINT = {\url{https://math.stackexchange.com/q/1422718}},
	URL = {\url{https://math.stackexchange.com/q/1422718}}
}

@article{hamkins,
        title={Pointwise definable models of set theory},
	volume=78, DOI={10.2178/jsl.7801090},
	number=1,
	journal={The Journal of Symbolic Logic},
	author={Hamkins, Joel David and Linetsky, David and Reitz, Jonas},
	year=2013,
	pages={139--156}
}

@article{hickman,
	title={{$\Lambda$}-Minimal Lattices},
	author={Hickman, John L.},
	year={1980},
	journal={Mathematical Logic Quarterly},
	volume={26},
	issue={10-12},
	pages={181-191},
	url={https://doi.org/10.1002/malq.19800261004},
	sn={0942-5616}
}

@MISC {stack_powerset,
	TITLE = {Hartogs number and the three power sets},
	AUTHOR = {Andrés E. Caicedo},
	HOWPUBLISHED = {MathOverflow},
	NOTE = {\url{https://mathoverflow.net/q/98381} (version: 2012-05-31)},
	EPRINT = {https://mathoverflow.net/q/98381},
	URL = {https://mathoverflow.net/q/98381}
}

@book{kunen1980,
   title = "Set theory: An Introduction to Independence Proofs",
   author = "Kunen, Kenneth",
   series = "Studies in Logic and the Foundations of Mathematics",
   publisher = "Elsevier Science",
   address = "Amsterdam, Lausanne, New York",
   url = "http://opac.inria.fr/record=b1117475",
   isbn = "0-444-86839-9",
   year = 1980
}

@ARTICLE{2023arXiv231009643B,
       author = {{Blass}, Andreas and {Kulshreshtha}, Dhruv},
        title = "{Cardinal Well-foundedness and Choice}",
      journal = {arXiv e-prints},
         year = 2023,
        month = oct,
          eid = {arXiv:2310.09643},
        pages = {arXiv:2310.09643},
          doi = {10.48550/arXiv.2310.09643},
archivePrefix = {arXiv},
       eprint = {2310.09643},
 primaryClass = {math.LO},
       adsurl = {https://ui.adsabs.harvard.edu/abs/2023arXiv231009643B}
}

\end{document}